\documentclass{article}

\usepackage{amssymb,amsmath,xcolor,graphicx,xspace,colortbl,textcomp, rotating} 
\usepackage{textcomp}
\usepackage{endnotes}  
\graphicspath{{multdeny_graphics/}{multdeny_tcache/}{multdeny_gcache/}}
\DeclareGraphicsExtensions{.pdf,.eps,.ps,.png,.jpg,.jpeg}

\title{Location and scale behaviour of the quantiles of a natural exponential family }

\author{Mauro Piccioni\thanks{Dipartimento di Matematica, Sapienza Universit\`{a} di Roma, 00185 Roma, Italia.
\texttt{piccioni@mat.uniroma1.it}},   Bartosz Ko\l odziejek \thanks{Wydzia{\l} Matematyki i Nauk Informacyjnych, Politechnika Warszawska, Warszawa, Polska. \texttt{b.kolodziejek@mini.pw.edu.pl}} and
G\'{e}rard Letac\thanks{Institut de Math\'ematiques de Toulouse, \'Equipe de Probabilit\'es et Statistique, Universit\'e  Paul Sabatier, 31062 Toulouse, France. \texttt{gerard.letac@math.univ-toulouse.fr}}}

=-0.5 in
\begin{document}
\maketitle 
\begin{abstract} Let $P_0$ be a probability on the real line  generating a natural exponential family  
$(P_t)_{t\in \mathbb {R}}$. Fix $\alpha$ in $ (0,1).$ We show that the property that $P_t((-\infty,t)) \leq \alpha \leq P_t((-\infty,t])$ for all $t$ implies that there exists a number $\mu_\alpha$ such that $P_0$ is the Gaussian distribution $N(\mu_{\alpha},1).$ In other terms, if for all $t$, $t$ is a quantile of $P_t$ associated to some threshold $\alpha\in (0,1)$, then the exponential family must be Gaussian. 
The case $\alpha=1/2$, \textit{i.e.} $t$ is always a median of $P_t,$ has been considered in Letac \textit{et al.} (2018). Analogously let $Q$ be a measure on $[0,\infty)$ generating a natural exponential family $(Q_{-t})_{t>0}$. We show that $Q_{-t}([0,t^{-1}))\leq \alpha \leq Q_{-t}([0,t^{-1}])$ for all $t>0$ implies that there exists a number $p=p_{\alpha}>0$ such that $Q(dx)\propto x^{p-1}dx,$ and thus $Q_{-t}$ has to be a gamma distribution with parameters $p$ and $t.$

\textsc{Keywords:}   Characterization of normal and gamma laws, one-dimensional exponential families, quantiles of a distribution, Deny equations.

\vspace{2mm}
\noindent
\textsc{MSC2010 classification:}  62E10, 60E05, 45E10.
\end{abstract}

\section{Introduction}

Let $P_0$ be a probability on the real line and  assume  its moment generating function
\begin{equation} \label{TLP}
M(t)=\int_{-\infty}^{+\infty}e^{tx}P_0(dx)
\end{equation}  
is finite for all real $t.$ Such a probability generates the natural exponential family  
\begin{equation}\label{NEP}
P_t(dx)=\frac{e^{tx}}{M(t)}P_0(dx),\  t\in \mathbb {R}.
\end{equation} 
For example, the Gaussian probability $P_0(dx)=(2\pi)^{-1/2}e^{-(x-m)^2/2}$, \textit{i.e.} $P_0=N(m,1)$  generates the natural exponential family
 $(P_t)=(N(m+t,1))$. In this case,  if $X_t \sim P_t$ for any $t \in \mathbb {R}$, then $X_t \sim X_0+t$, in other words $(P_t)$ is a location family generated by $P_0$ with location parameter $t$. It is well known and easily verified that the property $ X_t \sim X_0+t$ forces a natural exponential family to be generated by $P_0=N(m,1)$ for some $m.$  A way to see this is  to compute the m.g.f. of $X_t$ and substitute $X_0+t$ to $X_t$, getting the equation
$$
M(t+s)=M(t)M(s)e^{ts}, t, s \in \mathbb {R}.
$$
Taking logs and deriving w.r.t. $t$ and $s$ we get that the cumulant generating function $k=\log M$ of $P_0$ satisfies $k''(u)=1$ for all $u \in \mathbb {R}$, from which $k(u)=-mu+u^2/2$, that is precisely the c.g.f. of $N(m,1)$. 

The following remark is quite natural: the assumption that $X_t\sim X_0+t$, for any $t \in \mathbb {R}$, means that the distribution function of $X_t-t$ is independent of $t$, and so the same is true for the quantile function. If we make the weaker assumption that for some fixed $\alpha \in (0,1)$  an $\alpha$--quantile of $X_t-t$ does not depend on $t$, does one obtains the same
characterization? In slightly simplified words, if  $X_t \sim P_t$, as defined in \eqref{NEP}, is such that  $\Pr(X_t\leq t+b)=\alpha$  for any $t \in \mathbb{R}$, does this imply that $P_0$ is $N(m,1)$ for some $m?$

A recent paper (Letac \textit{et al} (2018)) gives the answer to this question for $\alpha=1/2$ (and $b=0$). Indeed it is  proved there that if $t$ is a median of $P_t$, for any $t \in \mathbb {R}$, then $P_0$ is the standard Gaussian $N(0,1)$. The first result of the present paper is the extension of this result for any $\alpha \in (0,1)$ (and an arbitrary $b$).

\vspace{4mm}\noindent \textbf{Theorem 1}. 
Let $P_0$ be a probability on the real line which generates the exponential family \eqref{NEP}. Let $b \in \mathbb {R}$ and suppose that $b+t$ is an $\alpha$-quantile of $P_t$, for $t \in \mathbb {R}$, that is
\begin{equation}\label{locationquant}
\int _{\left ( -\infty  ,b+t\right )}e^{tx}P_0(dx) \leq\alpha M(t) \leq \int _{\left ( -\infty  ,b+t\right ]}e^{tx}P_0(dx) ,\,\,\, t \in \mathbb {R} .
\end{equation}
Then $P_0=N(m^*,1)$, where $m^*=b-\Phi^{-1}(\alpha)$, $\Phi$ being the standard Gaussian distribution function; moreover $P_t=N(m^*+t,1)$.

\medskip

We describe now our second result: Let $Q$ be a Radon measure on the positive real line such that its Laplace transform
\begin{equation} \label{TLP2}
L(t)=\int_{[0,+\infty)}e^{-ty}Q(dy)
\end{equation}  is finite for all $t>0.$ 
Such a measure generates the exponential family  
\begin{equation}\label{NEP2}
Q_{-t}(dy)=\frac{e^{-ty}}{L(t)}Q(dy),\  t>0
\end{equation} 
(in principle this could be smaller than the natural exponential family generated by $Q$, but this will turn out to be impossible). For example, the  measure 
\begin{equation}\label{radon}
Q^{p}(dy)=\frac {1}{\Gamma (p)}y^{p-1}dy,
\end{equation}
defined for $p>0$ generates the natural exponential family $Q^p_{-t}=\text {Ga}(p,t)$, with $t>0$, where $\text {Ga}(p,t) $ is the gamma law with parameters $p$ and $t.$ Now it is immediately verified that if $Y_t \sim Q^p_{-t}$ then $Y_t \sim Y_1/t$,  that is $(Q^p_{-t})$ is a scale family generated by $Q^p_{-1}=\text {Ga}(p,1)$, with scale parameter $t^{-1}$. It is relatively easy to verify that this property forces $Q$ to be of the form \eqref{radon}. However the argument in the scale case  is slightly more involved than in the location  case and  we prefer to give the statement as a proposition.

\vspace{4mm}\noindent \textbf{Proposition 1}. Suppose that $(Q_{-t})_{t>0}$ is the natural exponential family defined in \eqref{NEP2}, for some  measure $Q$ on the non-negative real line. With $Y_t \sim Q_{-t}$, assume that $Y_t \sim \frac {Y_1}{t}$ for any $t>0$. Then, up to a multiplicative constant, $Q=Q^p$ defined by \eqref{radon}, for some $p>0$.

\vspace{4mm}\noindent\textbf{Proof of Proposition 1}. Compute the Laplace transform of $Y_t$ in the point $st$, where $s,t>0$. Then using the assumption $Y_t \sim \frac {Y_1}{t}$ one arrives at
$$
\frac {L(t+ts)}{L(t)}=\frac {L(1+s)}{L(1)}
$$
Defining $c(t)=\log L(t)$, for $t>0$ and deriving w.r.t. $t$ and $s$ this implies
$$
uc''(u)+c'(u)=0,
$$
where $u=t+ts>0$. Integrating twice one arrives at $c(u)=-p\log u+\ell$, with $p>0$ and an arbitrary $\ell \in \mathbb {R}$, from which $L(u)= \frac {e^{\ell}}{u^p}$, the Laplace transform of $e^{\ell}Q^p$. $\square$

It is worth to notice that the statement of the previous proposition and the analogous result for location families are special cases of the general results obtained by Ferguson (1962), that characterize \emph {general} exponential families which are location and scale families.

Now the assumption $Y_t \sim \frac {Y_1}{t}$ for any $t>0$, is equivalent to say that the distribution function of $tY_t$ is independent of $t$, and so the same is true for the quantile function. If we make the weaker assumption that, for some fixed  $\alpha \in (0,1)$,  an $\alpha$-quantile of $tY_t$ does not depend on $t$, is it enough to obtain the characterization stated in Proposition 1? In slightly simplified words, if $Y_t \sim Q_{-t}$ as defined in \eqref{NEP2}, is such that $\Pr(Y_t\leq a/t)=\alpha$  for all $t>0$, for some $a>0$, does this still imply that $Q$ is proportional to $Q^{p}$ for  some $p>0?$ Our second result gives a positive answer to this conjecture.

\vspace{4mm}\noindent \textbf{Theorem 2}. 
Let $Q$ be a  Radon measure on the non-negative real line which generates the exponential family \eqref{NEP2}. Let $a>0$ and suppose that $a/t$ is an  $\alpha$-quantile of $Q_{-t}$, for $t>0$, that is
\begin{equation}\label{scalequant}
\int _{\left ( 0 ,a/t\right )}e^{-ty}Q(dy) \leq\alpha L(t) \leq \int _{\left ( 0 ,a/t\right ]}e^{-ty}Q(dy) ,\,\,\, t>0 .
\end{equation}
 Then $Q$ is proportional to $Q^{p^*}$, where $p^*=p^*(\alpha)$ is the unique solution in $p>0$ of the equation $E_p(a)=\alpha$, $E_p$ being the distribution function of $\text {Ga}(p,1)$. In addition $Q_{-t}=\text {Ga}(p^*,t/a)$.

\vspace{4mm}
It is convenient to comment on the  existence and the uniqueness of $p^*$. The family $(\text {Ga}(p,1), p>0)$ is a convolution semigroup of laws supported by $(0,\infty).$ Hence, for any fixed $a>0$ the function  $p\mapsto E_p(a)$ is strictly decreasing in $p$ and is continuous. From the Markov inequality and the fact that the expectation for  $\text {Ga}(p,1)$ is $p$ we have $1- E_p(a)\leq p/a$ and this implies  $\lim_{p\downarrow 0}E_p(a)=1.$ The limit of $E_p(a)$ as $p \rightarrow \infty$ is zero from the law of large numbers.

\medskip

The proofs of Theorem 1 and 2 are given in the next section. These proofs  deduce from \eqref{locationquant} and \eqref{scalequant} two convolution equations in additive and multiplicative forms, respectively. The solutions to these equations have been investigated by Deny (1960). The result for additive convolutions is reported in the final section of  Deny (1960). The result for multiplicative convolutions can obtained with a passage to the additive convolution form by taking logarithms. In the next proposition we report both of them explicitly. 

\vspace{4mm}\noindent \textbf{Proposition 2}
\begin{itemize}
\item[1)] Suppose $H$ is a probability density defined on the whole real line, 
and consider the equation
\begin{equation}\label{Deny1}
f(t)=\int_{-\infty}^{+\infty} H(t-x)f(x)dx, t \in \mathbb {R},
\end{equation}
where $f$ is a locally integrable, non-negative function. Then $f$ is necessarily a linear combination, with non-negative coefficients, of a constant function with an exponential function of the form $e^{-s^*x}$, where $s^* \neq 0$ is a solution of the following equation in the real unknown $s$
\begin{equation}\label{moment}
\int_{-\infty}^{+\infty}e^{sx} H(x)dx=1.
\end{equation}
If there is no solution of this form then $f$ is necessarily constant.
\item[2)] Suppose $K$ is a probability density on the positive real line and consider the equation
\begin{equation}\label{Deny2}
g(t)=\int_{0}^{+\infty} K(\frac {t}{y})\, g(y)\, \frac {dy}{y}, \ \ t>0,
\end{equation}
where $g$ is a locally integrable and non-negative function on $(0,\infty)$. Then $g(t)$ is necessarily a linear combination, with non-negative coefficients, of the function $t^{-1}$ with a power function of the form $t^{-1-u*}$, where $u^* \neq 0$ is a solution of the following equation in the real unknown $u$
\begin{equation}\label{moment2}
\int_{0}^{+\infty}y^u K(y)dy=1.
\end{equation}
If there is no solution of this form then $g(t)=c/t$, where $c\geq 0$.
\end{itemize}

Both the equations \eqref{moment} and \eqref{moment2} have at most one non zero solution in $s$ and $u$, respectively. This follows by convexity of the logarithm of the functions appearing at the l.h.s. of these equations.

\section{Proofs}

\vspace{4mm}\noindent\textbf{Proof of Theorem 1}. 

Let us prove the theorem with $b=0$. Then we will adjust the solution to take into account an arbitrary value of $b$.
First we prove that $P_0$ is absolutely continuous. Take $ -A \leq s <t \leq A$, for some constant $A >0$ and compute

$$
P_0((s ,t)) =\int _{(s ,t)}e^{ -tx}e^{tx}P_0(dx) \leq e^{A^{2}}\int _{\left (s ,t\right )}e^{tx}P_0(dx)
$$
$$
\leq e^{A^{2}}\left (\int _{\left ( -\infty  ,t\right)}e^{tx}P_0(dx)-\int _{\left ( -\infty  ,s\right ]}e^{sx}P_0(dx) +\int _{\left ( -\infty  ,s\right ]}\left (e^{sx} -e^{tx}\right )P_0(dx)\right ).
$$
Using \eqref{locationquant} and the inequality $\left \vert e^{u} -e^{v}\right \vert  \leq \left \vert u -v\right \vert e^{w}$, for $\vert u\vert  ,\vert v\vert  \leq w$, this is bounded by
$$
e^{A^{2}}\left (\alpha \left (M(t) -M(s)\right ) +\vert t -s\vert \int _{\mathbb {R}}\vert x\vert e^{A\vert x\vert }P_0(dx)\right ) \leq c_{A}\vert t -s\vert 
$$
since $M$, being analytic, is locally Lipschitz, and the integral at the r.h.s. is finite by the existence of the m.g.f. of $P_0$ on the whole real line.

So we can always assume that $P_0$ has  a density $p_0$. Setting $\alpha=\frac {C}{1+C}$, with $C>0$, the quantile relation \eqref{locationquant} leads to
\begin{equation}\label{quanteq}
\int _{ -\infty }^{t}e^{tx}p_0(x)dx =C\int _{t}^{ +\infty }e^{tx}p_0(x)dx .
\end{equation}
Deriving w.r.t. $t$  both sides and multiplying by $e^{-t^{2}}$ one gets
\begin{equation}\label{anon}
p_0(t) +e^{ -t^{2}}\int _{ -\infty }^{t}xe^{tx}p_0(x)dx = -Cp_0(t) +Ce^{ -t^{2}}\int _{t}^{ +\infty }xe^{tx}p_0(x)dx
\end{equation}
Introduce the function defined by 
\begin{equation}\label{ABSC}\text {abs}_{C}(x) = -Cx1_{\left \{x <0\right \}} +x1_{\left \{x >0\right \}} .
\end{equation}
Multiply  both sides of  \eqref{quanteq} by $te^{-t^2}$ and subtract from \eqref{anon}. We obtain
\begin{equation}\label{unaeq}
p_0(t) =\frac{1}{1 +C}\int _{ -\infty }^{ +\infty }\text{abs}_{C}(t -x)e^{t(x -t)}p_0(x)dx .
\end{equation}
As expected, a solution to the equation \eqref{unaeq} is given by $\varphi (t -m^*)$, where $\varphi$ is the standard Gaussian density function, and 
$m^* = -\Phi ^{ -1}(\alpha)$. Next set
\begin{equation}\label{substi}
p_0(x) =\varphi (x -m)f(x),
\end{equation}
with $m=m^*$. We aim to prove that $f(x)$ has to be constant to solve the equation \eqref{unaeq}, with the substitution \eqref{substi}. Rewriting the equation for $f$, one gets
$$
f(t)e^{mt -t^{2}/2} =\frac{1}{1 +C}\int _{ -\infty }^{ +\infty }\text{abs}_{C}(t -x)e^{t(x -t) +mx -x^{2}/2}f(x)dx
$$
which is equivalent to
\begin{equation}\label{conveq}
f(t) =\frac{e^{\frac{m^{2}}{2}}}{1 +C}\int _{ -\infty }^{ +\infty }\text{abs}_{C}(t -x)e^{ -\frac{(t-x+m)^{2}}{2}}f(x)dx
\end{equation}
which has the form \eqref{Deny1} with
\begin{equation}\label{additivekern}
H(x) =\frac{e^{\frac{m^{2}}{2}}}{1 +C}\text {abs}_{C}(x)e^{ -\frac{(x +m)^{2}}{2}}.
\end{equation}
 The moment generating function of $H$ can be exactly computed
\begin{equation}\label{anon3}
\int _{ -\infty }^{ +\infty }e^{sx}H(x)dx =1 +\sqrt{2\pi }e^{(s-m)^2/2}(s-m)\left (\Phi (s-m)-\alpha \right ) .
\end{equation}
This is clearly equal to $1$ only if $s=0$ (hence $H$ is a density) and if  $s=m$. We apply Proposition 2 a) to the equation \eqref{conveq}. When $m=0$, that is if $\alpha=\frac {1}{2}$, the r.h.s. of \eqref{anon3}  is equal to $1$ only in $0$, hence the only non-negative non trivial solutions of the convolution equation \eqref{Deny1} with kernel $H$ given by \eqref{additivekern} are the positive constants. This yields immediately that $p_0(x)=\varphi (x)$, as desired. In the case $\alpha \neq \frac {1}{2}$ the solutions $f(x)$ are linear combinations with non-negative coefficients of the constant $1$ and the function $e^{-mx}$. Coming back to $p_0(x)=f(x)\varphi(x-m)$, this gives density solutions for $p_0$ which are mixtures of $N(m,1)$ with $N(0,1)$. But only the first one has the distribution function at $0$ equal to $\alpha \neq 1/2$, therefore a positive component from $N(0,1)$ is forbidden. This proves that $p_0(x)$ has  to be $\varphi(x-m)$.

Finally, to deal with an arbitrary value for $b$, define $\tau_{-b}=x-b.$  Now observe that if $P_t$ has $\alpha$-quantile $b+t$ then $P^*_t=P_t \circ \tau_{-b}^{-1}$ has $\alpha$-quantile $t$ and it is still a natural exponential family, for $t \in \mathbb {R}$. So $P^*_t=N(-\Phi^{-1}(\alpha),t)$ and $P_t=N(-\Phi^{-1}(\alpha)+b+t,1)$, ending the proof of Theorem 1. $\square$

\vspace{4mm}\noindent\textbf{Proof of Theorem 2}. 

 First we prove the result with $a=1$. Assume the relation \eqref{scalequant} and let $t \rightarrow 0+$. Then $L(t)$ increases to $Q(\mathbb {R}^{ +})$ by the monotone convergence theorem. Suppose $Q(\mathbb {R}^{ +})$ is finite: then, for any $\varepsilon  >0$, the l.h.s. of \eqref{scalequant} can be made larger than $Q(\mathbb {R}^{ +}) -\varepsilon $ in the following way. First choose $K$ in such a way that $Q(\left (0 ,K\right ]) >Q(\mathbb {R}^{ +}) -\frac{\varepsilon }{2}$; then choose $t <K^{ -1}$ and small enough to guarantee 
\begin{equation}
\int _{\left (0 ,t^{ -1}\right )}e^{ -tx}Q(dx) \geq \int _{\left (0 ,K\right ]}e^{ -tx}Q(dx) >Q(\mathbb {R}^{ +}) -\varepsilon  .
\end{equation}
Since $0<\alpha <1$, the first inequality in \eqref{scalequant} becomes absurd for $t$ and $\varepsilon$ sufficiently small. This implies that $Q(\mathbb {R}^{ +}) = +\infty $, hence the natural parameter space of the natural exponential family $(Q_{s})$ coincides with the negative reals.

Next we prove that $Q$ is absolutely continuous. Take $0 <s <t < +\infty $ and compute
\begin{eqnarray*}
Q((t^{ -1} ,s^{ -1})) &=&\int _{(t^{ -1} ,s^{ -1})}e^{sx}e^{ -sx}Q(dx)\leq e\int _{(t^{ -1} ,s^{ -1})}e^{ -sx}Q(dx)\\
&=&e\left (\int _{\left (0 ,s^{ -1}\right )}e^{ -sx}Q(dx) -\int _{\left (0 ,t^{ -1}\right ]}e^{ -tx}Q(dx) +\int _{\left (0 ,t^{ -1}\right ]}\left (e^{ -tx} -e^{ -sx}\right )Q(dx)\right ) .
\end{eqnarray*}
By \eqref{scalequant} the difference between the first two integrals at the r.h.s. is bounded by $\alpha \left (L(s) -L(t)\right )$, whereas the remaining integral is non positive. Again since $L$ is analytic in the positive real line it is locally Lipschitz and 
this proves the absolute continuity of $Q$, that is $Q(dx) =q(x)dx$, with $q$ non-negative and locally integrable.

Now we can write \eqref{scalequant} in the form of an equality, setting again $\alpha  =\frac{C}{1 +C}$, namely
\begin{equation}\label{nome}
\int _{0}^{t^{ -1}}e^{ -ty}q(y)dy =C\int _{t^{ -1}}^{ +\infty }e^{ -ty}q(y)dy , t>0.
\end{equation}
Deriving both sides w.r.t. $t$, one gets 
$$
\frac{1 +C}{t^{2}}e^{ -1}q(t^{ -1}) =C\int _{t ^{-1}}^{ +\infty }ye^{-ty}q(y)dy -\int _{0}^{t^{ -1}}ye^{ -ty}q(y)dy.
$$
Adding the l.h.s. of \eqref{nome} and subtracting the r.h.s., both multiplied by $t^{-1}$, to the r.h.s. of the above equality, we get for any $t>0$
\begin{equation} \label{anon5}
q\left (t^{ -1}\right ) =\frac{et^{2}}{1 +C}\left \{\int _{0}^{t^{ -1}}(t^{ -1} -y)e^{ -yt}q(y)dy +C\int _{t^{ -1}}^{ +\infty }(y -t^{ -1})e^{ -yt}q(y)dy\right \} .
\end{equation}
With the help of the function $\text{abs}_{C}$ defined  in \eqref{ABSC} , equality \eqref{anon5} is rewritten as
\begin{equation}\label{horrible}
q(t^{ -1}) =\frac{et}{(1 +C)}\int_0^{\infty} \text{abs}_{C}\left (1-ty\right )e^{ -ty}q(y)dy ,t >0.
\end{equation} 
Next, for any $p >0$, define $q_{p}(x) =\frac{1}{\Gamma (p)}x^{p -1}1_{\left (0 , +\infty \right )}(x)$. Recall from the introduction that for $p=p^*(\alpha)$ one has
\begin{equation}\label{forgotten}
\int_0^1 y^{p^*-1}e^{-y}dy=C\int_1^{+\infty}y^{p^*-1}e^{-y}dy
\end{equation} 
Now multiply both sides of \eqref{horrible} by $t^{p^* -2}$ and change the variable of integration at the r.h.s. to be $z =y^{ -1}$. One gets
\begin{equation}\label{encore}
t^{p^* -2}q(t^{ -1}) =\frac{et^{p^* -1}}{1 +C}\int_0^{\infty} \text{abs}_{C}\left (1-tz^{ -1} \right )e^{ -t/z}h\left (z^{ -1}\right )\frac{dz}{z^{2}} .
\end{equation}
Defining the l.h.s. of the above equality to be $g(t)$, one has
\begin{equation}\label{define}
q(t)=g(t^{-1})t^{p^*-2},
\end{equation}
and turns the equation \eqref{encore} into an equation of the form \eqref{Deny2} in $g$ with
\begin{equation}\label{chei}
K(y) =\frac{e}{1 +C}\text{abs}_{C}(1-y)e^{ -y}y^{p^* -1}1_{\left (0 , +\infty \right )}(y).
\end{equation}
The Mellin transform of $K$ can be easily computed
\begin{equation}\label{mellin}
\int_0^{\infty} y^u K(y)dy=1+e\Gamma(p^*+u)(p^*+u-1)\{\alpha-E_{p^*+u}(1)\}.
\end{equation}
Now observe that the quantity inside the brackets of \eqref{mellin} at the r. h. s. of \eqref{mellin} is always increasing in $u$; moreover it is equal to $0$ for $u=0$,
due to \eqref{forgotten}. Hence for any value of $C>0$ the function $K$ is always a density. When $p^*=1$ (equivalently, $C=e-1$, or $\alpha=1-e^{-1}$), $u=0$ is the unique global minimum point of the r. h. s. of \eqref{mellin}. Then, by Proposition 2 b), the only non negative non trivial solutions to the equation \eqref{Deny2}, 
with $K$ given by \eqref{chei}, have necessarily the form $g_0(t)=c_0t^{-1}$, with $c_0>0$. Thus $q(y)=c_0q_{p^*}(y)=c_0' y^{p^*-1}$.
Moreover, for $p^* \neq 1$ (equivalently, $C \neq e-1$, or $\alpha \neq 1-e^{-1}$) the value $u=1-p^* \neq 0$
makes the expression at the r. h. s. of \eqref{mellin} equal to $1$, too. As a consequence $g_1(t)=t^{p^*-2}$ is also a solution of the multiplicative convolution equation with $K$ given by \eqref{chei}. Applying again Proposition 2  b) all the non negative non trivial solutions are linear combinations 
$$
g_0(t)+c_1g_1(t)=c_0t^{-1}+c_1t^{p^*-2},
$$ 
with $c_0, c_1\geq 0$ not both zero. Substituting in \eqref{define} we get that the solutions to \eqref{anon5} have the form
$$
q(y)=c_0'y^{p*-1}+c_1.
$$ 
But for $c_1>0$ the condition \eqref{nome} is violated for $t=1$. Indeed, in this case the difference between the l.h.s. and the r.h.s. of \eqref{nome} is equal to $c_1(1-\frac {1+C}{e})$ and this is different from $0$  as soon as $C \neq e-1$. So again $q(y)=c_0'y^{p^*-1}$, as desired.

Finally, to deal with an arbitrary value of $a>0$, first define $\sigma_{a^{-1}}$ to be the multiplication by $a^{-1}$. Now observe that, if $Q_{-t}$ has $\alpha$-quantile $at^{-1}$ then $Q^*_{-t}=Q_{-t} \circ \sigma_{a^{-1}}^{-1}$ has $\alpha$-quantile $t^{-1}$ and it is still a natural exponential family, for $t>0$. So $Q^*_{-t}=\text {Ga}(p^*,t)$ and $Q_{-t}=\text {Ga}(p^*,at^{-1})$, ending the proof of Theorem 2.  $\square$

\section{References}
\vspace{4mm}\noindent \textsc{Deny, J.} (1960). Sur l'\'equation de convolution $\mu=\mu \ast \sigma$. 
\textit{S\'eminaire Brelot-Choquet-Deny (Th\'eorie du Potentiel)} \textbf{4}e ann\'ee, 1959-60,  
Expos\'e num\'ero 5.

\vspace{4mm}\noindent \textsc{Ferguson, T.S.} (1962), Location and scale parameters in exponential families of distributions, \textit{Ann. Math. Statist.} \textbf{33}, pp. 986-1001.

\vspace{4mm} \noindent \textsc{Letac, G., Mattner, L., Piccioni, M.} (2018), The median of an exponential family and the normal law, \textit{Statist. Prob. Lett.} \textbf{133}, pp. 38-41.

\end{document}